\input xy
\xyoption{all}

\input amssym.def
\input amssym.tex

\hfuzz=10pt

\magnification=\magstep1
 
 at10pt

\def\fitt{\hbox{\rm Fitt}}
\def\reg{\hbox{\rm reg}}

\def\ann{\hbox{\rm ann}}

\def\indeg{\hbox{\rm indeg}}

\def\proj{\hbox{\rm Proj}}

\def\sym{\hbox{\rm Sym}}

\def\C{{\cal C}}

\def\Z{{\cal Z}}
\def\R{{\cal R}}

\def\ip{{\frak p}}
\def\Ip{{\frak P}}

\def\im{{\frak m}}

\def\su{{\cal S}}
\def\H{{\cal H}}

\def\lra{{\longrightarrow}}

\def\a{{\alpha}}

\ \bigskip

\centerline{\bf Implicitization using approximation complexes}
\bigskip\bigskip\bigskip

\centerline{Marc Chardin}\smallskip
\centerline{Institut de Math{\'e}matiques de Jussieu}
\centerline{4 place Jussieu, Paris, France}
\medskip
\centerline{chardin@math.jussieu.fr}
\bigskip\bigskip\bigskip

{\bf 1. Introduction}\bigskip

We present in this short account a method to compute the image of a
rational map from ${\bf P}^{n-1}$ to ${\bf P}^{n}$, under suitable
hypotheses on the base locus and on the image. 

The formalism we use is due to Jean-Pierre Jouanolou, who gave a
course on this approach in the University of Strasbourg during the
academic year 2000--2001. In his joint article with Laurent Bus{\'e} [BJ],
this formalism is explained in details and applications to the 
implicitization problem are given. 

The idea of using a matrix of syzygies for the implicitization problem
goes back to the work of Sederberg and Chen [SC] and was at the origin
of several important contributions to this approach (see for instance
[Co], [CSC], [CGZ], [D'] and the articles on this subject in the volume of
the 2002 conference on Algebraic Geometry and Geometric Modeling
[AGGM02]). 

Most of this note is dedicated to presenting the method, the geometric
ideas behind it and the tools from commutative algebra that are
needed. Some references to classical textbooks are given for the
concepts and theorems we use for the presentation. In the last
section, we give the most advanced results we know related to this
approach. We illustrate this technique on an example that we carry out
in details all along the article.

References are given to the publication that fits best our 
statments. They may not be the first place where a similar result
appeared --for instance, many results were first proved for $n=2$ or
$n=3$--. \bigskip

{\bf 2. General setting}\bigskip

Given 

$$
\xymatrix{
\phi\ : {\bf P}^{n-1}\ar@{-->}[r]&{\bf P}^{n},\\}
$$
a rational map defined by $f:=(f_{0},\ldots
,f_{n})$,  $f_{i}\in R:=k[X_{1},\ldots ,X_{n}]$ homogeneous of degree
$d\geq 1$, such that the closure of its image is a hypersurface $\H$, the
goal is to compute the equation $H$ of this hypersurface.
\bigskip

We let :\medskip

$\bullet$ $I:=(f_{0},\ldots ,f_{n})\subset R$ be the ideal generated by the
$f_{i}$'s, \medskip

$\bullet$  $X:=\proj (R/I)\subset {\bf P}^{n-1}$ be the subscheme 
defined by $I$.\medskip
\bigskip

{\it A specific example :} We will illustrate in this article the
different steps and constructions on an example, taken from
[BCD, Example 3.2] that we will call example E : 

$$
\xymatrix{
\phi\ : {\bf P}^{2}\ar@{-->}[r]&{\bf P}^{3},\\}
$$

given by $f:=(ac^{2},b^{2}(a+c),ab(a+c),bc(a+c))$ with $R:={\bf
Q}[a,b,c]$. The ring of the target will be $R':={\bf
Q}[x,y,z,t]$.\bigskip

{\bf 3. The algebro-geometric intuition}
\bigskip

If $\Gamma_{0}\subset {\bf P}^{n-1}\times {\bf P}^{n}$ is the graph of
$\phi :({\bf P}^{n-1}-X)\lra {\bf P}^{n}$ and $\Gamma$ the Zariski
closure of $\Gamma_{0}$, one has :

$$
\H =\overline{\pi (\Gamma_{0})}=\pi (\Gamma )
$$
where $\pi :{\bf P}^{n-1}\times{\bf P}^{n}\lra {\bf P}^{n}$ is the
projection, and the bar denotes the Zariski closure (or equivalently
the closure for the usual topology in the case $k={\bf C}$). 

The first equality directly follows from the definition of $\H$, and
the second from the fact that $\pi$ is a projective morphism (so that
the image of a variety is a variety).\bigskip

On the algebraic side [Ha, II \S 7], one has
$$
\Gamma =\proj (\R_{I})
$$
with $\R_{I}:=R\oplus I\oplus I^{2}\oplus\cdots$ and the embedding
$\Gamma \subset {\bf P}^{n-1}\times{\bf P}^{n}$ corresponds to the
natural graded map :
$$
\eqalign{
S:=R[T_{0},\ldots ,T_{n}]&\buildrel{s}\over{\lra} \R_{I}\cr
 T_{i}\quad &\longmapsto f_{i}\in I=(\R_{I})_{1}\cr}.
$$

If $\Ip:=\ker (s)$, $\Ip_{1}$ (the degree 1 part of $\Ip$) is the
module of syzygies of the $f_{i}$'s :  
$$
a_{0}T_{0}+\cdots +a_{n}T_{n}\in \Ip_{1}\ \Longleftrightarrow\
a_{0}f_{0}+\cdots +a_{n}f_{n}=0.
$$

Setting  $\su_{I}:=\sym_{R}(I)$ and $V:=\proj (\su_{I})$, we have
natural onto maps 
$$
S\lra S/(\Ip_{1})\quad\quad\hbox{and}\quad\quad \su_{I}\simeq
S/(\Ip_{1})\lra S/\Ip\simeq \R_{I} 
$$
which corresponds to the embeddings
$$
\Gamma \subseteq V\subset {\bf P}^{n-1}\times{\bf
  P}^{n}.
$$

As $\R_{I}$ is the bigraded domain defining $\Gamma$, the projection
$\pi (\Gamma )$ is defined by the graded domain $\R_{I}\cap
k[T_{0},\ldots ,T_{n}]$. We have assumed that $\pi (\Gamma )$ is the
hypersurface $H=0$, so that this may be rewritten :
$$
(H)=\Ip\cap k[T_{0},\ldots ,T_{n}].
$$

{\sl In example E, with $S:=R[x,y,z,t]={\bf Q}[a,b,c,x,y,z,t]$ :
$$
\Ip =(ay-bz,at-cz,bt-cy,act-b(a+c)x)+(bx(z+t)-at^{2})+(xy(z+t)-zt^{2})
$$
where we have separated the (minimal) generators of degrees 1, 2 and 3
for simplicity. Of course it follows that $H=xy(z+t)-zt^{2}$.
Also, by definition, 
$$
\su_{I}\simeq S/(\Ip_{1})=S/(ay-bz,at-cz,bt-cy,act-b(a+c)x).
$$
}

The fact that $\R_{I}$ and $\su_{I}$, as well as the canonical map
$\su_{I}\lra\R_{I}$, do not depend on generators
of $I$ are useful to prove the following :\medskip

{\bf Theorem.} {\it $\Gamma =V$ if $X$ is locally a complete
intersection.}\medskip

{\sl In example E, the saturation of $I$ is the complete
  intersection ideal $(ac^{2},b(a+c))$ --because
  $b(a+c).(a,b,c)\subset I$ so that $I\subseteq
  (ac^{2},b(a+c))\subseteq I^{sat}$, and  $(ac^{2},b(a+c))$ is
  saturated--. Therefore $X$ is locally a   complete intersection (is
  it even globally a complete   intersection).} 
\bigskip

More refined criteria exist to insure that $\Gamma =V$, but we will
stick here to this one. This is partly justified by the following
result :\medskip

{\bf Proposition.} {\it If $\dim X=0$,  $\Gamma =V$ if and only if $X$ is
locally a complete intersection.}\medskip

The theorem above explains the key role of syzygies in computing
$H$ : they are equations of definition of $\Gamma$ when $X$ is locally
a complete intersection.\medskip

A more algebraic way to state the theorem is the following
:\medskip

{\bf Theorem.} {\it The prime ideal $\Ip$ is the saturation of the
  ideal generated by   its elements of degree 1 in the $T_{i}$'s (the
  syzygies) if $X$ is   locally a complete intersection.}\medskip  

{\sl In example E, $\Ip_{2}\subset (\Ip_{1}):(a,b,c)$ and $\Ip
  =(\Ip_{1}):(a,b,c)^{2}$.}\medskip

Nevertheless, as it is clear from this other formulation of the
theorem, one should not forget that even if $\Gamma =V$, it need
not be the case that $\R_{I}=\su_{I}$. In fact, the equality
$\R_{I}=\su_{I}$ may only hold in trivial cases in our context,
because $H$ is a minimal generator of $\Ip$. The difference between 
these algebras (which is the torsion part of $\su_{I}$, because
$\R_{I}$ is a domain) is a key point when one uses the syzygies to
compute $H$. This is very much similar to the fact that a homogeneous
ideal defining a variety in the projective space need not be
saturated.\medskip

The way the method proceeds is somehow parallel to determinantal methods for
computing resultants : it uses graded pieces of a resolution of
$\su_{I}$ to compute $\pi (V)$.\medskip

The connection between the elimination theory viewpoint, which looks
at $H$ as the generator of $\Ip\cap k[T_{0},\ldots ,T_{n}]$, and the
determinantal approach that computes $H$ from graded pieces of a 
resolution of $\su_{I}$ is shown by the following :\medskip

{\bf Proposition.} [BJ, 5.1] {\it Assume that $\Gamma =V$ and let
  $\eta$ be such that $H^{0}_{\im}(\su_{I})_{\mu}=0$ for all
  $\mu\geq\eta$. Then,
$$
\ann_{ k[T_{0},\ldots ,T_{n}]}(\su_{I}^{\eta})=\Ip\cap k[T_{0},\ldots ,T_{n}].
$$
}

Here $\im :=(X_{1},\ldots ,X_{n})$ and $H^{0}_{\im}(M):=\{ m\in M\ |\
\exists \ell ,\ X_{i}^{\ell}m=0 \ \forall i\}$.  
The graded pieces of $\su_{I}$ will be described below, and we will
provide estimates for $\eta$ satifying the vanishing condition. Notice that
$H^{0}_{\im}(\su_{I})$ is the torsion part of  $\su_{I}$ when  $\Gamma
=V$.\medskip

{\bf Remark.} The choices on gradings are one of the delicate points
in this approach. For instance, the hypothesis
$H^{0}_{\im}(\su_{I})_{\mu}=0$ is equivalent to 
$$
H^{0}_{\im}(\sym^{j}_{R}(I))_{\mu +dj}=0,\quad \forall j,
$$
if we adopt the natural grading of $\sym^{j}_{R}(I)$ making the
canonical map $\sym^{j}_{R}(I)\lra I^{j}\subset R$ a homogeneous map
of degree zero.
\medskip

A candidate for a resolution of $\su_{I}$ is the  $\Z$-complex
introduced and studied by Herzog, Simis and Vasconcelos. We will
decribe this complex in the next section.\bigskip  

{\bf 4. The tools from commutative algebra}\bigskip

{\it The saturation of an ideal.} --- An ideal $I$ in a polynomial
ring $R:=k[X_{0},\ldots ,X_{n}]$ is {\it saturated} (or, more precisely
$\im$-saturated) if $I:\im =I$, where $\im :=(X_{0},\ldots
,X_{n})$. In other words, $I$ is saturated if : $X_{i}f\in I,\ \forall
i\ \Rightarrow f\in I$.  

The ideal $I^{sat}:=\bigcup_{j}(I:\im^{j})$ is saturated, it is the
smallest saturated ideal containing $I$ and is called the saturation
of $I$. 


Another way of seeing the saturation of an ideal, that directly
extends to modules, is given by the remark that :
$$
I^{sat}=I+H^{0}_{\im}(R/I)
$$
that one can also write $R/I^{sat}=(R/I)/H^{0}_{\im}(R/I)$. 

The saturation of a module $M$ will be $M/H^{0}_{\im}(M)$. As usual,
one should be careful about the fact that the saturation of an ideal
$I$ then corresponds to saturating the module $R/I$ and not the 
ideal considered as a module over the ring.

Seeing the saturation operation in relation with the left exact
functor $H^{0}_{\im}(\hbox{---})$ naturally leads to the consideration
of the derived functors $H^{i}_{\im}(\hbox{---})$, and to the
cohomological approach of algebraic geometry. 

There is a one-to-one correspondance between the subschemes of a
  projective space ${\bf   P}^{n}_{k}$ and the {\it saturated}
  homogeneous ideals of the polynomial ring $R:=k[X_{0},\ldots
  ,X_{n}]$. To see this notice that, by definition, two subschemes of
  ${\bf   P}^{n}_{k}$ are the same if they coincide on all the affine charts
  $X_{i}=1$. If $\phi_{i}$ is the specialization homomorphism
  $X_{i}\mapsto 1$ then  the homogeneization of $\phi_{i}(I)$ is the
  ideal $I_{(i)}:=\bigcup_{j}(I:(X_{i}^{j}))$. It follows that $I$ and
  $J$ defines the same schemes if and only if $I_{(i)}=J_{(i)}$ for
  all $i$, which is easily seen to be equivalent to the equaltiy of
  their saturation as $I^{sat}=\bigcap_{i}I_{(i)}$.

When considering multigraded ideals, with respect to set of variables
that are generating ideals $\im_{1},\ldots ,\im_{t}$
(these ideals are never maximal unless $t=1$ and $k$ is a field), the
operations of saturation with respect to the different ideals
naturally appears. The subschemes of the corresponding product of
projective spaces corresponds one-to-one to ideals that are saturated
with respect to all the ideals $\im_{i}$, or equivalently with respect
to the product of these ideals.\medskip 

{\it The ring of sections.} [Ha, II \S 5, Ei \S A4.1] --- If
$R:=k[X_{0},\ldots ,X_{n}]$ and $B:=R/I$ is the quotient of $R$ by the
homogeneous ideal $I$, an interesting object to consider is : 
$$
\Gamma B:=\ker \big( \bigoplus_{i}B_{X_{i}}\lra
  \bigoplus_{i<j}B_{X_{i}X_{j}}\big) ,
$$
where $B_{(f)}:=\{ {{x}\over{f^{j}}}\ |\ x\in B,\ j\in {\bf N}\}$ and
  the maps are the evident ones up to a sign chosen so that $(1,\ldots
  ,1)$ maps to $0$. One has a natural isomorphism $B_{X_{i}}\simeq
  B/(X_{i}-1)$ and $\Gamma B$ should be interpreted as the
  applications that are defined on each affine chart $X_{i}=1$ and
  matches on the intersection of any two of these charts. Notice that it
  is clear from the definition that replacing $I$ by its saturation do
  not affect $\Gamma B$.

In a sheaf theoretic language, one has
$$
\Gamma B=\bigoplus_{\mu\in {\bf Z}}H^{0}({\bf P}^{n}_{k},{\cal
  O}_{X}(\mu )),
$$
with $X:=\proj (B)$, and the natural grading of $\Gamma B$ coincides
with the grading of the section ring on the right hand side.

These considerations extends to modules along the same lines. Also, the
map we used to define $\Gamma B$ fits into a complex, called the \v
Cech complex,
$$
\xymatrix@R=8pt{
0\ar[r]&B\ar^(.4){\phi}[r]\ar@{=}[d]&\bigoplus_{i}B_{X_{i}}\ar^(.4){\psi}[r]\ar@{=}[d]&
\bigoplus_{i<j}B_{X_{i}X_{j}}\ar[r]\ar@{=}[d]&\cdots\ar[r]&B_{X_{0}\cdots
  X_{n}}\ar[r]\ar@{=}[d]&0\\
&{\cal C}^{0}_{\im}(B)&{\cal C}^{1}_{\im}(B)&{\cal C}^{2}_{\im}(B)&&
{\cal C}^{n+1}_{\im}(B)&\\}.
$$

One has $H^{0}_{\im}(B)=\ker (\phi )$ and $\Gamma B=\ker (\psi )$. It
is a standard fact that $H^{i}_{\im}(B)$ is isomorphic to the $i$-th
cohomology module of this complex. This in particular gives an exact
sequence :  
$$
\xymatrix{
0\ar[r]&H^{0}_{\im}(B)\ar[r]&B\ar[r]&\Gamma B\ar[r]&H^{1}_{\im}(B)\ar[r]&0\\},
$$
which splits into two parts the difference between the homogeneous
quotient $B$ and the more geometric notion of the section ring
attached to $X:=\proj (R/I)\subseteq {\bf P}^{n}_{k}$. 

Notice that, in the case $k$ is a field and $X$ is of dimension zero,
$\Gamma B_{\mu}=H^{0}({\bf P}^{n}_{k},{\cal   O}_{X}(\mu ))$ is a
$k$-vector space of dimension the degree of $X$ for any $\mu$. In
particular, when $\dim X=0$, $\Gamma B$ is not finitely
generated. In any dimension, it can be shown that $\Gamma (R/I)$ is
finitely generated if anf only if $I$ has no associated prime $\ip$
such that $\proj (R/\ip)$ is of dimension zero ({\it i.e.} $\dim
(R/\ip )=1$).\medskip

{\it Castelnuovo-Mumford regularity.} [Ei, \S 20.5] --- The
Castelnuovo-Mumford 
regularity is an invariant that measures the algebraic complexity
of a graded ideal or module over a polynomial ring $R:=k[X_{0},\ldots 
,X_{n}]$. The two most standard definitions are given either in terms
of a minimal finite free $R$-resolution of the module (this
resolution exists by Hilbert's theorem on syzygies) or in terms of the
vanishing of the cohomology modules defined above (using a theorem of
Serre [Ha, III 5.2] to show that this makes sense). \medskip 

{\bf Theorem-Definition.} {\it Let $b_{i}(M)$ be the maximal degree of a
minimal $i$-th syzygy of $M$ and $a_{i}(M):=\inf\{\mu\in {\bf Z}\ |\
H^{i}_{\im}(M)_{\nu}=0,\ \forall \nu >\mu \}$, then 
$$
\reg (M)=\max_{i}\{ a_{i}(M)+i\} =\max_{i}\{ b_{i}(M)-i\} .
$$
}

Notice that if $M=R/I$, minimal $0^{\rm th}$ syzygies of $M$ are
minimal generators of $M$ (namely, the element $1$), minimal
$1^{\rm st}$
syzygies of $M$ are minimal generators of $I$, and $2^{\rm nd}$
syzygies of $M$ are syzygies between the chosen (minimal) generators of
$I$. If one looks at $I$ as a module, these modules are the same up to
a shift in the labeling, except $0^{\rm th}$ module for $R/I$, and one
has $\reg (I)=\reg (R/I)+1$.  

The existence of different interpretations of the regularity is a key
to many results on this invariant. It is for instance immediate from the
cohomological definition that $\reg (I^{sat})\leq \reg (I)$, but this
is not easy to see on the definition in terms of syzygies.
Also, when $\dim X=0$ ($X:=\proj (R/I)$, as above), it easily follows
from the cohomological definition and the fact that $H^{i}_{\im}(M)=0$
for $i>\dim M$ (Grothendieck's vanishing theorem) that $\reg (I)$ is
the smallest integer $\mu$ such that :\smallskip   

(1) $I_{\mu}=(I_{X})_{\mu}$ (recall that $I_{X}$ is the saturation of $I$),\smallskip

(2) $\dim (R/I_{X})_{\mu -1}=\deg (X)$.\medskip

In case $X$ is a set of simple points, condition (2) says that passing
through the $\deg (X)$ different points of $X$ impose linearly
independant conditions on polynomials of degree $\mu -1$. An
elementary account on regularity in this context is given in \S 4 of
[AGGM02, D. Cox. Curves, surfaces, and syzygies, 131--150].  

The fact that $\reg (I)$ bounds the degrees of the syzygies of $I$
shows the naturality of considering this invariant in the
implicitization problem using the syzygy matrix.

On the computational side, the degrees of generators of a
Gr{\"o}bner basis of the ideal for the degree-reverse-lex order, under a
quite weak conditions on the coordinates, is bounded by $\reg
(I)$. This is another way of understanding the regularity as a measure of
the complexity of the ideal.

\medskip

{\it The Koszul complex.} [Ei, \S 17] --- Let $x=(x_{1},\ldots ,x_{r})$ be a $r$-tuple of
elements in a ring $A$. The (homological) Koszul complex
$K_{\bullet}(x;A)$ is the complex with modules
$K_{p}(x;A):=\bigwedge^{p}A^{r}\simeq A^{{p}\choose{r}}$ and maps
$d_{p}:K_{p}(x;A)\lra K_{p-1}(x;A)$ defined by :  
$$
g.e_{i_{1}}\wedge\cdots \wedge e_{i_{p}}\longmapsto g.\sum_{j=1}^{p}(-1)^{j+1}x_{i_{j}}
e_{i_{1}}\wedge\cdots \widehat{\wedge e_{i_{j}}}\cdots \wedge e_{i_{p}}.
$$

We set $Z_{p}(x;A):=\ker (d_{p})$ and $H_{p}(x;A):=Z_{p}(x;A)/{\rm
  im}(d_{p+1})$.\bigskip


{\it The $\Z$-complex.} [Va, Ch. 3] --- We consider $f_{i}\in
R\subset S$ as elements of $S$ and the two complexes
$K_{\bullet}(f;S)$ and $K_{\bullet}(T;S)$ where $T:=(T_{0},\ldots
,T_{n})$. These complexes have the same modules
$K_{p}=\bigwedge^{p}S^{n+1}\simeq S^{{p}\choose{n+1}}$ and 
differentials $d^{f}_{\bullet}$ and $d^{T}_{\bullet}$\medskip

$\bullet$ It directly follows from the definitions that
$d^{f}_{p-1}\circ d^{T}_{p}+d^{T}_{p-1}\circ d^{f}_{p}=0$, so that 
$d^{T}_{p}(Z_{p}(f;S))\subset Z_{p-1}(f;S)$. The complex
$\Z_{\bullet}:=(Z_{\bullet}(f;S),d^{T}_{\bullet})$ is the called
$\Z$-complex associated to the $f_{i}$'s.\medskip

$\bullet$ Notice that $Z_{p}(f;S)=S\otimes_{R}Z_{p}(f;R)$ and \medskip

--- $Z_{0}(f;R)=R$, \medskip

--- $Z_{1}(f;R)={\rm Syz}_{R}(f_{0},\ldots ,f_{n})$,\medskip

--- the map $d_{1}^{T}:S\otimes_{R}{\rm Syz}_{R}(f_{0},\ldots
,f_{n})\lra S$ is defined by :
$$
(a_{0},\ldots ,a_{n})\longmapsto a_{0}T_{0}+\cdots +a_{n}T_{n}.
$$

The following result shows the intrinsic nature of the homology of the
$\Z$-complex, it is a key point in proving results on its
acyclicity.\medskip 

{\bf Theorem.} {\it $H_{0}(\Z_{\bullet})\simeq \su_{I}$ and the
homology modules $H_{i}(\Z_{\bullet})$ are $\su_{I}$-modules that only
depends on $I\subset R$, up to isomorphism.}\medskip

$\bullet$ We let $R':=k[T_{0},\ldots ,T_{n}]$ and look at graded pieces :
$$
\Z_{\bullet}^{\mu}\ :\ \cdots\lra 
R'\otimes_{k}Z_{2}(f;R)_{\mu}\buildrel{d_{2}^{T}}\over{\lra}
R'\otimes_{k}Z_{1}(f;R)_{\mu}\buildrel{d_{1}^{T}}\over{\lra} 
R'\otimes_{k}Z_{0}(f;R)_{\mu}\lra 0
$$ 
where $Z_{p}(f;R)_{\mu}$ is the part of  $Z_{p}(f;R)$ consisting of
elements of the form $\sum a_{i_{1}\cdots i_{p}}e_{i_{1}}\wedge\cdots
\wedge e_{i_{p}}$ with the $a_{i_{1}\cdots i_{p}}$ all of the same
degree $\mu$.\medskip

{\bf Nota Bene.} {\it This is not the usual convention for the grading of these
  modules, however we choosed it here for simplicity. The usual
  grading (used for instance in {\rm [BC]} or {\rm [Ch]}) makes the Koszul maps
  homogeneous of degree 0, so they ask $a_{i_{1}\cdots i_{p}}$ to be
  homogeneous of degree $\mu -pd$ in place of being of degree $\mu$.}\medskip

We will denote the cokernel of the last map by $\su_{I}^{\mu}$. 

\bigskip 


{\it Determinants of complexes.} [No \S 3.6, GKZ App. A] --- Let $A$
be a commutative domain, for simplicity.\medskip

If $A^{n}\buildrel{\alpha}\over{\lra}A^{n}$ is $A$-linear we can
define $\det (\alpha )\in A$.\medskip

If we have a complex $C_{\bullet}$ with three terms:

$$
\xymatrix@R=3pt@C=35pt{
  &A^{m}\ar[r]&A^{n}\\
  &\bigoplus& \\
A^{m}\ar[r]\ar^{\beta}[uur]&A^{n}\ar^{\alpha}[uur]\\}
$$
such that $\det (\beta )\not= 0$, we set $\det (C_{\bullet}):={{\det (\alpha
    )}\over{\det (\beta )}}$. In fact $\det (C_{\bullet})$ is
independent of the decomposition of $C_{1}$ as a direct sum
$A^{m}\oplus A^{n}$. \medskip

More generally a bounded complex $C_{\bullet}$ of free $A$-modules such that
${\rm Frac A}\otimes_{A}C_{\bullet}$ is exact may always be decomposed in the
following way

$$
\xymatrix@R=3pt@C=35pt{
&\cdots&A^{q}\ar[r]  &A^{p}\ar[r] &A^{m}\ar[r]&A^{n}\\
& &\bigoplus &\bigoplus&\bigoplus& \\
\cdots&A^{q}\ar[r]\ar^{\delta}[uur]&A^{p}\ar[r]\ar^{\gamma}[uur]&A^{m}\ar[r]\ar^{\beta}[uur]&A^{n}\ar^{\alpha}[uur]\\}
$$
with $\alpha, \beta ,\ldots$ having non-zero determinants. Then $\det
(C_{\bullet}):={{\det (\alpha ).\det (\gamma )\cdots}\over{\det     (\beta
    ).\det (\delta )\cdots }}$ is independent of the decomposition.\medskip

Performing the decomposition of a given complex is easy : decompose first
$C_{1}$ into $A^{n}\bigoplus A^{m}$ so that $\det (\a )\not= 0$ (this
amounts to choose a non zero maximal minor of the map $C_{1}\lra
C_{0}=A^{n}$), and apply the procedure recursively to the complex : 
$$
\xymatrix@R=3pt@C=35pt{
\cdots\ar[r]&C_{3}\ar[r]&C_{2}\ar[r]
&A^{m}\ar[r]&0\\}
$$
\bigskip

{\it Fitting ideals.} [No \S 3.1, Ei, \S 20.2] --- If $A$ is a ring and $M$ is a module
represented as the cokernel of a 
map $\psi :A^{m}\lra A^{n}$, the ideal generated by minors of size
$n-i$ of $\psi$ only depends on $M$ and $i$, this ideal is called the
$i$-th Fitting ideal of the $A$-module $M$. One of the most important of
these ideals associated to the $A$-module $M$ is the $0$-th Fitting ideal
({\it i.e.} the one generated by the maximal minors of $\psi$),
denoted by $\fitt^{0}_{A}(M)$.

\bigskip
{\bf 5. The method and main results}\bigskip

We assume hereafter that $\pi (\Gamma )$ is of codimension 1 in ${\bf
  P}^{n}_{k}$ defined by the equation $H=0$ and denote by $\delta$ the
  degree of the map $\pi$ from $\Gamma$ onto its image. 
\bigskip

If $J$ is a $R'$-ideal, we will denote by $[J]$ the gcd of the
elements in $J$. It represents the component of codimension one of the
scheme defined by $J$ (its divisorial component) because $R'$ is
factorial.\medskip

With these notations, one has :\medskip

{\bf Proposition 1.} [BJ, 5.2] {\it If $X=\emptyset$, $\Z_{\bullet}$ is acyclic and 
$$
[\fitt^{0}_{R'}(\su_{I}^{\mu})]=\det (\Z_{\bullet}^{\mu})=H^{\delta},
$$
for every $\mu \geq (n-1)(d-1)$.}\medskip

The identities above are identities of principal ideals in $R'$,
therefore it corresponds to an equality of elements of $R'$ up to
units. Recall that $[\fitt^{0}_{R'}(\su_{I}^{\mu})]$ is the gcd of
maximal minors of the map  $d_{1}^{T}:(a_{0},\ldots ,a_{n})\longmapsto
a_{0}T_{0}+\cdots +a_{n}T_{n}$ from the syzygies of degree $\mu$
(each $a_{i}$ is of degree $\mu$) seen as a vector space over $k$ to
$R'\otimes_{k}R_{\mu}$. The entries of this matrix are therefore
linear forms in the $T_{i}$'s with coefficients in $k$.

\medskip

This proposition shows that the determinant of this graded part of
$\Z_{\bullet}$ actually computes the divisor $\pi_{*}(\Gamma )=\delta
.\pi (\Gamma )$ obtained as direct image of the cycle $\Gamma$ (see
[Fu, \S 1.4] for the defintion of the direct image $\pi_{*}(\Gamma )$ of the
cycle $\Gamma$). 

\medskip

In the case $X$ is of dimension zero, the situation is slightly more
complicated :\medskip

{\bf Proposition 2.} [BJ, 5.7, 5.10 \&\  BC 4.1] {\it If $\dim X=0$,\medskip

{\rm (i)} The following are equivalent :\smallskip

$\quad$ {\rm (a)} $X$ is locally defined by at most $n$ equations,\smallskip

$\quad$ {\rm (b)} $\Z_{\bullet}$ is acyclic,\smallskip

$\quad$ {\rm (c)} $\Z_{\bullet}^{\mu}$ is acyclic for $\mu\gg 0$.\medskip

{\rm (ii)} If $\Z_{\bullet}$ is acyclic, then
$$
[\fitt^{0}_{R'}(\su_{I}^{\mu})]=\det (\Z_{\bullet}^{\mu})=H^{\delta}G,
$$
for every $\mu \geq (n-1)(d-1)-\varepsilon_{X}$, where
$1\leq \varepsilon_{X}\leq d$ is the minimal degree of a hypersurface
containing $X$ and $G\not= 0$ is a homogeneous polynomial which is a
unit if and only if $X$ is locally a complete intersection.}\bigskip 

{\bf Remark 3.} In fact $[\fitt^{0}_{R'}(\su_{I}^{\mu})]=\det
(\Z_{\bullet}^{\mu})=\pi_{*}V$ for $\mu \geq
(n-1)(d-1)-\varepsilon_{X}$, and the degree of $G$ is the sum of
numbers measuring how far $X$ is from a complete intersection at each
point of $X$.\medskip

{\bf Remark 4.} It is very fast to compute the ideal $I_{X}$ with a
 dedicated computer algebra system (like Macaulay 2, Singular or
Cocoa), and {\it a fortiori} to compute $\varepsilon_{X}$ which is the
 smallest degree of an element in $I_{X}$. Moreover the following
 result actually implies a good bound on the complexity of this
 task.\medskip

{\bf Proposition 5.} [Ch, 3.3] {\it If $J\subset R$ is a homogeneous ideal generated
in degree at most $d$ with $\dim (R/J)=1$ and $J'$ its saturation (in
other words, the defining ideal of the zero-dimensionnal scheme
$X:=\proj (R/J)$), then 
$$
\reg (J)\leq n(d-1)+1\quad{\rm and}\quad \reg (J')\leq (n-1)(d-1)+1.
$$
}

{\sl In example E, $\reg (I)=\reg (I_{X})=4$, while the general bound
above gives $\reg (I)\leq 7$ and $\reg (I^{sat})\leq 5$.
A minimal free $R$-resolution of $I$ gives a resolution of $Z_{1}$ :
$$
\xymatrix{
0\ar[r]&R[-2]\ar^(.35){\left[ \matrix{c\cr -b\cr -a\cr
      0\cr}\right]}[r]&R[-1]^{3}\oplus R[-2]\ar^(.5){\left[ \matrix{
0 &0 &0 &-b(a+c)\cr
a &0 &c&0\cr 
-b&-c&0 &0\cr 
0 &a &-b &ac\cr}
\right]}[rrrrr]&&&&&Z_{1}\subset R^{4}=S_{1}\\} 
$$
and we have seen that $I_{X}=(ac^{2},b(a+c))$, so that $\indeg
(I_{X})=2$ and therefore $(n-1)(d-1)-\varepsilon_{X}=2\times (3-1)-2=2$.
The syzygies of degree 2 are of the form :
$$
\ell_{1}(ay-bz)+\ell_{2}(at-cz)+\ell_{3}(cy-bt)+\lambda_{4}(act-b(a+c)x)
$$
with $\ell_{i}\in R_{1}$  and $\lambda_{4}\in R_{0}=k$. Notice
that they are not linearly independant, and that the relation (unique
in this degree) is given by the second syzygy : 
$$
c(ay-bz)-b(at-cz)-a(cy-bt)=0.
$$
We may for instance choose as generators of syzygies of degree 2 the 9
syzygies, $s_{1}$ to $s_{9}$ :
$$
a(ay-bz),b(ay-bz),a(at-cz),b(at-cz),c(at-cz),a(cy-bt),b(cy-bt),c(cy-bt),act-b(a+c)x
$$
which gives the $6\times 9$ matrix of linear forms (elements of
$R'_{1}$) for the matrix of $d_{1}^{T}$ in degree 2 (recall that
$T=(x,y,z,t)$ with the notations of the example) :
$$
\matrix{
\matrix{\hphantom{ab}}
\quad\ &
\matrix{
s_{1}&s_{2}&s_{3}&s_{4}&s_{5}&s_{6}&s_{7}&s_{8}&s_{9}\cr
\hphantom{-z}&\hphantom{-z}&\hphantom{-z}&\hphantom{-z}&\hphantom{-z}&
\hphantom{-z}&\hphantom{-z}&\hphantom{-z}&\hphantom{-z}\cr}\cr
\matrix{
a^{2}\cr
ab\cr
ac\cr
b^{2}\cr
bc\cr
c^{2}\cr}
\quad\,&
\left[
\matrix{
 y& 0& t& 0& 0& 0& 0& 0& 0\cr
-z& y& 0& t& 0&-t& 0& 0&-x\cr
 0& 0&-z& 0& t& y& 0& 0& t\cr
 0& -z& 0& 0& 0& 0&-t& 0& 0\cr
 0& 0& 0&-z& 0& 0& y&-t&-x\cr
 0& 0& 0& 0&-z& 0& 0& y& 0\cr}
\right] .\cr}
$$

Now, $Z_{2}$ has a free $R$-resolution of the form:
$$
\xymatrix{
0\ar[r]& R[-4]\ar[r]&R[-1]\oplus R[-3]^{3}\ar[r]&Z_{2}\subset 
\bigwedge^{2}S_{1}= R^{6}\\}. 
$$
In particular, $Z_{2}$ has one minimal generator of degree 1 and no
minimal generator of degree 2. The element of degree 1
$$
\Sigma := c\; y\wedge z -b\; z\wedge t-a\; y\wedge t \in \bigwedge^{2}S_{1}
$$
satisfies $d_{2}^{f}(\Sigma )=b(a+c)[c(ay-bz)-b(at-cz)-a(cy-bt)]=0$. Therefore
$(Z_{2})_{1}=\Sigma k$ and $(Z_{2})_{2}=\Sigma R_{1}$. 
We have $d_{2}^{T}(\Sigma )=c(z\otimes y-y\otimes z)-b(z\otimes t-t\otimes z)-a(
t\otimes y-y\otimes t)\in S\otimes_{R}\bigwedge^{1}S_{1}$,
that we may rewrite
$$
d_{2}^{T}(\Sigma )= -t\otimes (ay-bz) +y\otimes(at-cz) +z\otimes (cy-bt).
$$
In degree 2, the matrix of $d_{2}^{T}:R'\otimes_{k}\Sigma R_{1}\lra
R'\otimes_{k}Z_{1}(f;R)_{2}$ on the bases $(a\Sigma, b\Sigma, c\Sigma )$
for the source and $(s_{1},\ldots ,s_{9})$ for the target is therefore
the transpose of
$$
\matrix{
\matrix{\hphantom{a\Sigma}}\quad\ &
\matrix{
s_{1}&s_{2}&s_{3}&s_{4}&s_{5}&s_{6}&s_{7}&s_{8}&s_{9}\cr
&&&&&&&&\cr}\cr
\matrix{
a\Sigma \cr
b\Sigma \cr
c\Sigma \cr}
\quad\,&
\left[
\matrix{
-t&0&\, y\,&0&0&z&0&0&0\cr
0&-t&0&y&0&0&z&0&0\cr
0&0&0&-t&\, y\,&\, t\,&\ 0\,&\, z\,&0\cr
}
\right] .\cr}
$$
We now choose a maximal non zero minor of this matrix, for instance
the minor $\Delta_{2}$ given by lines 3, 4 and 5 of the matrix of
$d_{2}^{T}$, and the minor $\Delta_{1}$ of the matrix of $d_{1}^{T}$
obtained by erasing columns 3,4 and 5. We get the formula : 

$$
H={{\Delta_{1}}\over{\Delta_{2}}}=
{{\left\vert 
\matrix{
 y& 0& 0& 0& 0& 0\cr
-z& y& t& 0& 0&-x\cr
 0& 0& y& 0& 0& t\cr
 0& -z&0&-t& 0& 0\cr
 0& 0& 0& y&-t&-x\cr
 0& 0& 0& 0& y& 0\cr}
\right\vert}
\over
{\left\vert 
\matrix{
y&0&0\cr
0&y&-t\cr
0&0&y\cr}
\right\vert}}
={{-y^{3}(xyz+xyt-t^{2}z)}\over{y^{3}}}.
$$
}

{\it Computations of the free $R$-resolutions of $Z_{1}$ and $Z_{2}$ were
done using the dedicated software {\rm Macaulay 2} by Dan Grayson and
Mike Stillman {\rm [M2]}. In the case $n=3$, this computation goes
very fast, even for pretty high degree $d$, and {\rm Macaulay 2}
performs degree truncations to speed up the computation, if
needed. The graded pieces that we need to know can also easily be
computed using linear algebra routines, as detailed in {\rm   [BC]}
and implemented in {\rm [Bu]}.}\medskip  

When the dimension of the base locus $X$ of the map $\phi$ increases,
the situation becomes harder to analyze. In dimension 1, the situation
is pretty well understood :\medskip

{\bf Proposition 6.} [Ch, 8.2, 8.3] {\it Assume that $\dim X=1$ and let
  $\C$ be the union of components of dimension 1 of $X$ (its ``unmixed
  part''). Then, 
\medskip

{\rm (i)} The following are equivalent :\smallskip

$\quad$ {\rm (a)} $X$ is locally defined by at most $n$ equations and $\C$
is defined on a dense open subset by at most $n-1$ equations,\smallskip

$\quad$ {\rm (b)} $\Z_{\bullet}^{\mu}$ is acyclic for $\mu\gg 0$.\medskip

{\rm (ii)} If $\Z_{\bullet}^{\mu}$ is acyclic for $\mu\gg 0$, the following
are equivalent :\smallskip

$\quad$ {\rm (a)}  $\Z_{\bullet}$ is acyclic,\smallskip

$\quad$ {\rm (b)} $\C$ is arithmetically Cohen-Macaulay,\smallskip

$\quad$ {\rm (b')} every section $f\in H^{0}(\C ,{\cal O}_{\C}(\mu ))$ is
the restriction to $\C$ of a polynomial function of degree $\mu$, for
every $\mu\in {\bf Z}$.\medskip

{\rm (iii)} If $\Z_{\bullet}^{\mu}$ is acyclic for $\mu\gg 0$ and  $H^{0}(\C
,{\cal O}_{\C}(\mu ))=0$ for all $\mu <-d$ ---for instance if $\C$ is
reduced--- then  $\Z_{\bullet}^{\mu}$ is acyclic for $\mu\geq
(n-1)(d-1)$. If further $X$ is defined by at most $n-1$ equations
locally on the support of $\C$, then 
$$
[\fitt^{0}_{R'}(\su_{I}^{\mu})]=\det (\Z_{\bullet}^{\mu})=H^{\delta}G,
$$
for every $\mu \geq (n-1)(d-1)$, where $G$ is a homogeneous polynomial
such that the support of $\pi (V)$ is the zero set of $GH$.}\medskip

Here also, more precisely, $\det (\Z_{\bullet}^{\mu})$ represents
the divisor $\pi_{*}V$.\medskip 

{\bf Remark 7.} It is perhaps true that  $\det (\Z_{\bullet}^{\mu})$ represents
the divisor $\pi_{*}V$ for $\mu \geq (n-1)(d-1)$ when
$\Z_{\bullet}^{\mu}$ is acyclic for $\mu\gg 0$ and  $H^{0}(\C ,{\cal
  O}_{\C}(\mu ))=0$ for all $\mu <-d$, but we needed the slightly
stronger hypothesis above to prove it in [Ch]. 

\bigskip\bigskip

{\bf References.}\medskip

{\it Articles.}\medskip

[AGGM02] R. Goldman, R. Krasaukas (Eds.). Topics in Algebraic Goemetry
and Geometric Modeling. {\it Contemporary Mathematics} {\bf 334} (2003).\medskip

[BJ] L. Bus{\'e}, J.-P. Jouanolou.  On the closed image of a
  rational map and the implicitization problem. {\it J. Algebra} {\bf
  265} (2003), 312--357.\medskip

[BC] L. Bus{\'e}, M. Chardin. Implicitizing rational hypersurfaces
  using approximation complexes. {\it J. Symbolic Computation (to
  appear).}\medskip

[BCD] L. Bus{\'e}, D. Cox, C. D'Andr{\'e}a.
Implicitization of surfaces in ${\bf P}^3$ in the presence of base
  points. {\it J. Algebra Appl.} {\bf 2} (2003), 189--214. \medskip 

[Ch] M. Chardin. Regularity of ideals and their powers. Preprint
{\bf 364} (Mars 2004), Institut de Math{\'e}matiques de Jussieu, Paris. \medskip 

[Co] D. Cox. Equations of parametric curves and surfaces via syzygies.
{\it Contemporary Mathematics} {\bf 286} (2001), 1--20.\medskip

[CSC] D. Cox, T. Sederberg, F. Chen.
The moving line ideal basis of planar rational curves.
{\it Comp. Aid. Geom. Des.} {\bf 15} (1998), 803--827.\medskip

[CGZ] D. Cox, R. Goldman, M. Zhang.
On the validity of implicitization by moving quadrics for rationnal
  surfaces with no base points.
{\it J. Symbolic Computation} {\bf 29} (2000), 419--440.\medskip

[D']  C. D'Andr{\'e}a. Resultants and moving surfaces.
{\it J. of Symbolic Computation} {\bf 31} (2001), 585--602.\medskip

[SC] T. Sederberg, F. Chen.
Implicitization using moving curves and surfaces.
{\it Proceedings of SIGGRAPH 95}, Addison Wesley, 1995,
301--308.\medskip

\medskip

{\it Textbooks.}\medskip

[Ei] D. Eisenbud. {\it Commutative algebra. With a view toward
  algebraic geometry.}  Graduate Texts in Mathematics {\bf
  150}. Springer-Verlag, New York, 1995. \medskip

[Fu] W. Fulton. {\it Intersection theory.} Second edition. Ergebnisse der
Mathematik und ihrer Grenzgebiete {\bf 3}. Springer-Verlag, Berlin, 
1998.  \medskip 

[GZK] I. M. Gel'fand, M. Kapranov, A. Zelevinsky. {\it Discriminants,
resultants, and multidimensional determinants.} Mathematics: Theory \&
Applications. Birkh{\"a}user Boston, Inc., Boston, MA, 1994. \medskip

[Ha] R. Hartshorne. {\it Algebraic geometry.} Graduate Texts in
Mathematics {\bf 5}2. Springer-Verlag, New York-Heidelberg, 1977. \medskip

[No] D. G. Northcott. {\it Finite free resolutions.} Cambridge Tracts
in Mathematics {\bf 7}1. Cambridge University Press, Cambridge-New
York-Melbourne, 1976. \medskip 

[Va] W. Vasconcelos. {\it The Arithmetic of Blowup Algebras}.
London Math. Soc. Lecture Note Ser. {\bf 195}. 
Cambridge University Press, 1994.\medskip

\medskip

{\it Softwares.} \medskip

[Bu] L. Bus{\'e}. {\it Algorithms for the implicitization using approximation
  complexes.} \break\noindent ({\sl
  http://www-sop.inria.fr/galaad/personnel/Laurent.Buse/program.html}).\medskip

[M2] D. Grayson, M. Stillman. {\it Macaulay 2.} 
({\sl http://www.math.uiuc.edu/Macaulay2/}).
\end